\UseRawInputEncoding
\documentclass[12pt]{article}
\usepackage{amssymb}
\usepackage{amsmath}
\begin{document}

\centerline{\Large \bf Existence and smoothness of the Navier-Stokes equations}

\centerline{\Large \bf and semigroups of linear operators}
\vskip 2cm
\centerline{\sc Yu.\,N.\,Kosovtsov}

\medskip
\centerline{Lviv, Ukraine}
\centerline{email: {\tt yunkosovtsov@gmail.com}} \vskip 1cm
\begin{abstract}
Based on Leray's formulation of the Navier-Stokes equations and the conditions of the exact linear representation of the nonlinear problem found in this paper, a compact explicit expression for the exact operator solution of the Navier-Stokes equations is given. It is shown that the introduced linear operator for Leray's equations is the generator of one-parameter contraction semigroup. This semigroup yields the existence of a unique and smooth classical solution of the associated Cauchy problem of Navier-Stokes equations in space $\mathbb{R}^3$ under smooth initial conditions.
\end{abstract}

\section{The Navier-Stokes Equations}

Incompressible flows of homogeneous fluids in space $\mathbb{R}^3$ are solutions of the system of equations \cite {Ladyzh}-\cite {Gilles}
\begin{equation}
\frac{\partial { v}}{\partial t} =\nu\Delta v-\sum_{j=1}^3 v_j \frac{\partial { v}}{\partial x_j} -\nabla p+f,
\label{N-Sv}
\end{equation}

\begin{equation}
 \mathrm{div} \; v \equiv \sum_{j=1}^3 \frac{\partial { v_j }}{\partial x_j}=0, \qquad  (x,t) \\ \in \\ \mathbb{R}^3 \times [0,\infty ),
\label{N-Svo}
\end{equation}

\begin{equation}
v\mid _{t=0}=v_0, \qquad  x \\ \in \\ \mathbb{R}^3 ,
\label{N-V0}
\end{equation}
where $v(x, t) \equiv (v_1, v_2, v_3 )$ is the fluid velocity, $p(x, t)$ is the scalar pressure, the gradient operator is $\nabla = (\frac{\partial}{\partial x_1},\frac{\partial}{\partial x_2},\frac{\partial}{\partial x_3})$ and the Laplace operator $\Delta\equiv\sum_{j=1}^3 \frac{\partial ^2}{\partial x_j^2}$, a kinematic constant viscosity $\nu \geq≥ 0$, initial conditions $v\mid _{t=0}=v_0(x) \equiv (v_{01}, v_{02}, v_{03} )$, $f=f(x,t)$ is given, externally applied force. Further we will consider Navier-Stokes equations in the absence of external forces, i.e. with $f=0$.

It is well known \cite {Majda},\cite {Foias} that we can find an equation $p = p(v)$ to eliminate the pressure from (\ref{N-Sv}). If we take the divergence at each side of the momentum equation (\ref{N-Sv}), which yields, using (\ref{N-Svo}) that
\begin{equation}
\Delta\emph{p}=-\sum_{i,j=1}^3 \frac{\partial { v_j}}{\partial x_i} \frac{\partial { v_i}}{\partial x_j}.
\label{eqPoisson}
\end{equation}

The equation (\ref{eqPoisson}) with respect to the function $p$ is the well known Poisson equation, whose solution in $\mathbb{R}^3$ \cite {Majda} has the form (for brevity, we will sometimes denote this solution as $p_{v}(x,t)\equiv p_{v}$)

\begin{equation}
 p(x,t)=\frac{1}{4 \pi} \int_{\mathbb{R}^3} d\xi \sum_{i,j=1}^3 \frac{\partial { v_j(t,\xi)}}{\partial \xi_i} \frac{\partial { v_i(t,\xi)}}{\partial \xi_j}\frac{1}{|x-\xi|}, \qquad ( \xi \\ \in \\ \mathbb{R}^3  ).
\label{Poissonsoln}
\end{equation}

Substituting the solution (\ref{Poissonsoln}) into (\ref{N-Sv}), we obtain (Leray’s Formulation of the Navier-Stokes Equation \cite {Majda}) \emph{non-linear system of integro-differential equations} containing \emph{only} unknown functions $v_i$.
\begin{equation}
\frac{\partial { v}}{\partial t} =\nu\Delta v-\sum_{j=1}^3 v_j \frac{\partial { v}}{\partial x_j} -\nabla  p_v, \qquad ( x \\ \in \\ \mathbb{R}^3 , t \geq 0 ).
\label{N-Svv}
\end{equation}

It is easy to see \cite {Majda} that the systems of the Navier-Stokes equations (\ref{N-Sv})-(\ref{N-V0}) are equivalent to the system of the equations (\ref{Poissonsoln}),(\ref{N-Svv}), provided that functions $v_0$ in initial conditions are smooth and
\begin{equation}
\sum_{i=1}^3 \frac{\partial { v_{0i}}}{\partial x_i}=0,  \qquad ( x \\ \in \\ \mathbb{R}^3  ).
\label{N-S20}
\end{equation}

Although the Leray's formulation of the Navier-Stokes equation has been known for a long time, it has not been widely used to analyze its solutions. The main theses of this can be traced to the following \cite {Majda}.

The equation (\ref{N-Svv}) is quadratically nonlinear and contains a nonlocal, quadratically nonlinear operator. These facts make the Navier-Stokes equation hard to study analytically.

Although (\ref{N-Svv}) constitutes a closed system for $v$, this formulation turns out to be not very useful for further analysis except by rather crude methods based on the energy principle. The main reason for this is that this formulation hides all the properties of vorticity stretching and interaction between rotation and deformation.

\section{Exact linear representation of a nonlinear problem}
We define for a multi-index $(\alpha_1, . . . , \alpha_1) \in \mathbb{N}^n $
\[D^\alpha = D_1^{\alpha_1} ... D_3^{\alpha_n},\; \mathrm{where}\; D_k = \frac{\partial }{\partial x_k},\;  \alpha := \sum_{i=1}^n\alpha_i. \]

Let us take an arbitrary function $u(x)$, $x \in  \mathbb{R}^n$ which is infinitely many times differentiable and decay sufficiently rapidly, i.e., $u \in C^\infty(\mathbb{R}^n)$, and

\[\lim_{\mid x\mid\rightarrow\infty}\mid x\mid^p D^\alpha u(x) = 0 \qquad \mathrm{for\; some} \;p \in \mathbb{N}\; \mathrm{and\; all }\; \alpha \in \mathbb{N}^n. \]

Then, using this function, we define the linear space of \emph{analytic} on all arguments functions $G(u):\equiv G(u,u_1,u_2,...)\in C^\infty$,
where $u_i:=D_k^i u(x)$ ($i,k \in \mathbb{N}$) and which decay sufficiently rapidly with $u$ and all $u_i$ provided $\mid x\mid\rightarrow\infty$.

Since we can also consider these functions $G(u(x))$ as functions of $x  \in  \mathbb{R}^n$ denoting $G(u(x))=\widetilde{G}(x)$. In this case, it is obvious that the function $\widetilde{G}(x) \in C^\infty(\mathbb{R}^n)$ and  decay sufficiently rapidly when $\mid x\mid\rightarrow\infty$.

When endowed with the usual inner product
\[<\widetilde{F}(x),\widetilde{G}(x)>:=\int_{\mathbb{R}^n}\widetilde{F}(x) \widetilde{G}(x) \; dx \]
the considered space becomes as a complete normed subspace $L^2_u(\mathbb{R}^n) $ of the Hilbert space $L^2(\mathbb{R}^n) $.

Let a differential operator $A$ with domain $D(A)$ on a Banach space  $L^2_u(\mathbb{R}^n) $ is of type

\begin{equation}
 A=\int_{\mathbb{R}^n} d \zeta F(u(\zeta))\frac{\delta}{\delta u(\zeta)}, \qquad  \zeta \\ \in \\ \mathbb{R}^n ,
\label{opA}
\end{equation}
where $u(x)$ is a function of class described above, a \emph{predetermined} $F(u(x)):\equiv F(u,u_1,u_2,...)\in  L^2_u(\mathbb{R}^n) $ with $u_i:=D_k^i u(x)$ and we can consider $F(u(x))$ as the  linear multiplication operator, $\frac{\delta}{\delta u(\zeta)}$ is the functional derivative.  We have specially transferred the symbol $d \zeta$ at the beginning of the formula to emphasize that the right operator is executed first,
and then the integration is performed.

Here by functional derivative we mean the linear mapping with the following property
\[\frac{\delta u(x)}{\delta u(\zeta)}=Dirac(x-\zeta),\]
the chain rule is also valid in this context.

Entered  operator is a linear operation, i.e., operator $A$ obeys the following properties:
\[A c G(u(x))= c A G(u(x)),  \qquad \mathrm{for\; any}\;G\in \\ L^2_u(\mathbb{R}^n)\;\mathrm{and\; any\; scalar}\; c\in \\ \mathbb{R}, \]
\[A [G(u(x))+H(u(x))]= A G(u(x))+A H(u(x)),\qquad \mathrm{for\; any}\; G,H\in \\ L^2_u(\mathbb{R}^n).\]

It is obvious that $A^ju\in L^2_u(\mathbb{R}^n)$ for any $j \in \mathbb{N}$ and operator $A$ is \emph{continuous} and therefore it is \emph{bounded} and \emph{closable} on the space $\L^2_u(\mathbb{R}^n)$.

Let an operator $A$ with some $F(u(x))$ is  the \emph{generator of strongly continuous semigroup} $(T(t))_{t\geq0}$ then (\cite {Engel}, p.50)
\begin{equation}
\frac{d}{dt}T(t)u(x) = T(t)Au(x) = AT(t)u(x) \qquad \mathrm{for \; all}\; t \geq 0.
\label{eq3}
\end{equation}

Note that if we denote
\begin{equation}
v(x,t) := T(t)u(x),
\label{solF}
\end{equation}
then from (\ref{eq3}) \emph{three} \emph{equivalent} equations follows.

First
\begin{equation}
\frac{d}{dt}v(x,t) = Av(x,t),
\label{lin}
\end{equation}
the \emph{linear} one with the \emph{functional derivative}  and expression (\ref{solF}) is its \emph{general} solution.

Second
\[\frac{d}{dt}v(x,t) = T(t)Au(x),\]
but because for considered operator \[Au(x)=F(u(x)),\] then as much as (see \cite {Engel}, p. 58) $T(t)$ is a semigroup of algebra homomorphisms on $D(A)$, i.e.,
\begin{equation}
T(t)(h \cdot g) =T(t)h \cdot T(t)g \qquad \mathrm{for}\; h, g \in  D(A)\; \mathrm{and}\; t \geq 0
\notag
\end{equation}
and, since we considering analytic functions $F$, we have \[T(t)Au(x)=F(T(t)u(x))\] and come to
\begin{equation}
\frac{d}{dt}v(x,t) = F(v(x,t),v_{x_1}(x,t),...,D^\alpha v(x,t),...),
\label{nonlin}
\end{equation}
i.e., to \emph{nonlinear PDE}. The expression (\ref{solF}) is the \emph{general} solution to (\ref{nonlin}) \emph{too}!

This result has parallels with the well known connection between \emph{nonlinear} ODEs and \emph{linear} PDEs.
There is a third equivalent equation
\[Av(x,t)=F(v(x,t),v_{x_1}(x,t),...,D^\alpha v(x,t),...),\]
the \emph{nonlinear} one with the \emph{functional derivative} and with the \emph{same} \emph{general} solution.

From a comparison of expressions (\ref{opA}) and (\ref{nonlin}), it is easy to obtain a one-to-one relationship between the linear operator $A$ and the nonlinear right-hand side of equation (\ref{nonlin}).

 Note that the strongly continuous semigroup $(T(t))_{t\geq0}$ generated by the operator $A$ of the type considered here with some $F(u(x))$ can be
given by \cite {Engel}, \cite {Ren}
\begin{equation}
T(t) = \exp\{t A\} =\sum_{n=0}^\infty \frac{t^n A^n}{n!}, \qquad t\geq 0,
\notag
\end{equation}
 where the power series converges for every $t \geq 0$ (\cite {Engel}, pp. 52, 81).

Practically without changing the above reasoning, we can generalize the results obtained to more complex equations, when \[F(v(t))=F(v(t,x),v_{x_1}(t,x),...,D^\alpha v(t,x),...,I_1v(t,x),I_2v(t,x)...),\] where $I_i$ are nonlinear \emph{integral} operators.

Moreover, it is possible to consider in a completely similar way systems of nonlinear equations

\begin{equation}
\frac{d}{dt}v_i(x,t) = F_i(v_1(x,t),...v_n(x,t)),\;\; i=1..n,
\label{sys}
\end{equation}
where all $F_i$ are functions of the type described above and if we introduce the following linear operator
\begin{equation}
 A=\int_{\mathbb{R}^n} d \zeta \sum_{i=1}^n F_i(u_1(\zeta),...,u_n(\zeta))\frac{\delta}{\delta u_i(\zeta)}, \qquad  \zeta \\ \in \\ \mathbb{R}^n,
\label{sysopA}
\end{equation}
where rapidly decreasing $u_i \in C^\infty(\mathbb{R}^n)$, and $F_i(u_1(x),...u_n(x))$ are in the space $\L^2_u(\mathbb{R}^n)$ built on  vector function $u(x) \equiv (u_1(x), ..., u_n(x) )$.

Again assuming that operator $A$ is the \emph{generator of strongly continuous semigroup} $(T_{sys}(t))_{t\geq0}$.
In this case the general solutions to the system (\ref{sys}) is
\begin{equation}
v_i(x,t) = T_{sys}(t)u_i(x),
\label{solsys}
\end{equation}
where as before
\begin{equation}
T_{sys}(t) = \exp\{t A\} =\sum_{n=0}^\infty \frac{t^n A^n}{n!}, \qquad t\geq 0.
\notag
\end{equation}
Thus we have obtained an important result that, under the above conditions, it is possible to reduce \emph{exactly} a \emph{nonlinear} problem to a \emph{linear} one, which allows us to use for solving nonlinear problems the accumulated baggage for linear problems\footnotemark. \footnotetext{Here we consider the case when the operator $A$ does not depend on $t$. It is possible to generalize the results obtained above for more general functions $F(t,x,u(x))$ (see the idea in \cite {Kos2}).}

\textbf{Comment}. So far as $T(t)\mid_{t=0}=1$ (identity operator) then from (\ref{solF}) or (\ref{solsys}) it follows that  function $u(x)$ is directly related to the initial conditions. However, we need to note one circumstance in previous reasoning. We use in definition of the operators $A$  the functional derivatives with respect to the variable $u$ and therefore cannot identify it with everyone fixed initial condition. That is to obtain the particular solution we have to substitute $u(x)=v(x,0)$ into the general solution
\begin{equation}
v(x,t) = \{T(t)u(x)\}\mid_{u(x)=v(x,0)}.
\label{psolF}
\end{equation}

\section{Solution of the Navier-Stokes problem}

Turning to the Leray's formulation of the Navier-Stokes equations (\ref{N-Svv}), we will consider the vector function $u(x) \equiv (u_1(x), u_2(x), u_3(x) )$, where arbitrary functions $u_i(x)$, $x \in  \mathbb{R}^3$ which, as we discussed in the previous section, are infinitely many times differentiable and decay sufficiently rapidly. In accordance with the trick described above let us introduce the \textit{linear} differential operator
\begin{align}
&{\bf A}=\int_{\mathbb{R}^3} d \zeta \,\sum_{i=1}^3 {\Large\{}\nu\Delta_{\zeta} u_i(\zeta)-\sum_{j=1}^3 u_j(\zeta) \frac{\partial { u_i(\zeta)}}{\partial \zeta_j}-\notag
\\\notag
\\&\frac{1}{4 \pi}\frac{\partial }{\partial \zeta_i}{\Large[} \int_{\mathbb{R}^3} d\xi \sum_{i,j=1}^3 \frac{\partial { u_j(\xi)}}{\partial \xi_i} \frac{\partial { u_i(\xi)}}{\partial \xi_j}\frac{1}{|\zeta-\xi|}{\Large]}{\Large\}}\frac{\delta}{\delta u_i(\zeta)}.
\label{OPERn}
\end{align}
In short notation
\begin{equation}
{\bf A}=\int_{\mathbb{R}^3} d \zeta \,\sum_{i=1}^3 {\Large\{}\nu\Delta_{\zeta} u_i(\zeta)-\sum_{j=1}^3 u_j(\zeta) \frac{\partial { u_i(\zeta)}}{\partial \zeta_j}-\nabla_{\zeta}\; p_u(\zeta) {\Large\}}\frac{\delta}{\delta u_i(\zeta)},
\label{OPERp}
\end{equation}
where
\[ p_u(\zeta)\equiv \frac{1}{4 \pi}{\Large[} \int_{\mathbb{R}^3} d\xi \sum_{i,j=1}^3 \frac{\partial { u_j(\xi)}}{\partial \xi_i} \frac{\partial { u_i(\xi)}}{\partial \xi_j}\frac{1}{|\zeta-\xi|}{\Large]}.\]

We restrict our attention to Navier-Stokes equations and correspondingly to operator ${\bf A}$ on all of $\mathbb{R}^3 $, thus avoiding the delicate discussion of boundary conditions.

As it follows from the results of the previous section we need to find out whether the operator is a generator of a semigroup.
Let us describe some properties of the operator ${\bf A}$ introduced above. First of all
\begin{equation}
{\bf A}\,u_i=\nu\Delta u_i-\sum_{j=1}^3 u_j \frac{\partial { u_i}}{\partial x_j}- \frac{1}{4 \pi}\frac{\partial }{\partial x_i} \int_{\mathbb{R}^3}d\xi \sum_{i,j=1}^3 \frac{\partial { u_j(\xi)}}{\partial \xi_i} \frac{\partial { u_i(\xi)}}{\partial \xi_j}\frac{1}{|x-\xi|}
\label{Ax}
\end{equation}
and it is obvious that
${\bf A}\;L^2_u(\mathbb{R}^3)\mapsto L^2_u(\mathbb{R}^3)$ and ${\bf A}$ is \emph{continuous} and therefore it is \emph{bounded} and \emph{closable} on the space $\L^2_u(\mathbb{R}^n)$ built on vector function $u(x) \equiv (u_1(x), u_2(x), u_3(x) )$.
We can rewrite the expressions (\ref{Ax}) in vector form for brevity
\begin{equation}
{\bf A}u =\nu\Delta u-\sum_{j=1}^3 u_j \frac{\partial { u_i}}{\partial x_j} -\nabla  p_u,
\label{Au}
\end{equation}
where
\begin{equation}
p_u(x)=\frac{1}{4 \pi}\int_{\mathbb{R}^3}d\xi \sum_{i,j=1}^3 \frac{\partial { u_j(\xi)}}{\partial \xi_i} \frac{\partial { u_i(\xi)}}{\partial \xi_j}\frac{1}{|x-\xi|}.
\notag
\end{equation}

Further, by analogy with the proof of the energy principle, we have in view of (\ref{N-S20}) and that $u$ is rapidly decreasing and hence $p_u$ rapidly decreasing at infinity too  one can obtain the well known result that inner product
\begin{equation}
<{\bf A}u,u>=-\nu\sum_{i,j=1}^3\int_{\mathbb{R}^3}{\Large(}\frac{\partial { c_i}}{\partial x_j}{\Large)}^2 dx\leq 0.
\label{Adis}
\end{equation}
Hence, the linear operator ${\bf A}$ on the Hilbert space $L^2_u(\mathbb{R}^3)$ is \emph{dissipative} \cite {Hille}-\cite {Phillips}.

By definition \cite {Dun}, the adjoint operator $A^*$ to the bounded operator $A$ on the Hilbert space satisfies the identity

\begin{equation}
<Au,h>=<u,A^*h>, \qquad \mathrm{for\; any} \; u,h\in L^2(\mathbb{R}^3)).
\notag
\end{equation}
In the case of the \emph{real} vector Hilbert space we have for considered operator ${\bf A}$, taking into account (\ref{Adis}), that
\begin{equation}
<u,{\bf A}^*u>=<{\bf A}^*u,u>=<{\bf A}u,u>\;\leq 0, 
\notag
\end{equation}
so adjoint operator ${\bf A}^*$ is \emph{dissipative} too.

Now we use one of the results of semigroup theory (\cite {Engel}, p. 84, 3.17 Corollary).

\textbf{Proposition 1}. \emph{Let} $(A,D(A))$ \emph{be a densely defined operator on a Banach
space} $X$. \emph{If both} $A$ \emph{and its adjoint} $A'$	\emph{ are dissipative, then the closure} $\overline{A}$ \emph{of}
$A$ \emph{generates a contraction semigroup on} $X$.

A contraction semigroup is a special case of strongly continuous semigroup $(T(t))_{t\geq0}$ \cite {Engel}, \cite {Ren}.
In each case, the semigroup is given by
\begin{equation}
T(t) = \exp\{tA\} =\sum_{n=0}^\infty \frac{t^n A^n}{n!}, \qquad t\geq 0,
\notag
\end{equation}
 where the power series converges for every $t \geq 0$.

So the \emph{general} solution to the Leray’s system (\ref{N-Svv}) is as follows
\begin{equation}
v(x,t) = \exp\{t{\bf A}\}u(x) =\sum_{n=0}^\infty \frac{t^n {\bf A}^n}{n!}u(x), \qquad t\geq 0,
\label{Ler}
\end{equation}
where $u(x)$ is arbitrary vector function.

And its particular solutions
\begin{equation}
v(x,t) = \{\exp\{t{\bf A}\}u(x)\}\mid_{u=v_0} =\{\sum_{n=0}^\infty \frac{t^n {\bf A}^n}{n!}u(x)\}\mid_{u(x)=v_0}, \qquad t\geq 0,
\label{LerNS}
\end{equation}
provided that $v_0:=v_0(x)$ is divergence-free vector field and the scalar pressure $p$ restored through known solution $v$ by (\ref{Poissonsoln}) leads us to the solution of the Navier-Stokes equations (\ref{N-Sv})-(\ref{N-V0}).

Because the operator ${\bf A}$ is the generator of a contraction semigroup, it follows that the semigroup yields solutions of the
associated abstract Cauchy problem (\cite {Engel}, p.145).
Furthemore, if $u(x)\in D({\bf A})$ the map $v(x,t) := T(t)u(x)$ is the \emph{unique classical} solution of the X-valued initial value problem (or abstract Cauchy problem).

Returning to the Navier-Stokes equations, we come to the conclusion that the operator solution of the system of equations in form (\ref{LerNS}) provided that initial condition $v_0=v(x,0)$ is divergence-free is the \emph{unique classical} solution to (\ref{N-Svv}).

That is, we have proved

\textbf{Proposition 2}. \emph{Let initial conditions} $v(x,0)$ \emph{be any smooth rapidly decreasing, divergence-free vector field}.
\emph{Take} $f(x, t)$ \emph{to be identically zero. Then there exist unique smooth rapidly decreasing functions} $p(x, t)$, $v_i(x, t)$
\emph{on} $\mathbb{R}^3\times  [0,\infty)$ \emph{that satisfy} (\ref{N-Sv})-(\ref{N-V0}).

I wish to thank Prof. Pierre Gilles Lemari\'{e}-Rieusset for his valuable comments on this article.

\end{document}